\documentclass{elsarticle}
\usepackage{graphicx}
\usepackage{float}
\usepackage{subfigure}
\usepackage{amsfonts}
\usepackage{amsmath}
\usepackage{amsthm}
\usepackage{amssymb}
\usepackage{bm}
\usepackage{color,todonotes}
\usepackage{multirow}
\usepackage[T1]{fontenc}
\usepackage[utf8]{inputenc}
\usepackage{comment}
\usepackage{mathtools}
\usepackage{lineno}

\newenvironment{system}%
	{\left\lbrace \begin{array}{@{} l @{} }}%
	{ \end{array}\right.}

\newcommand{\bs}{\boldsymbol}
\newcommand{\de}{\,\mathrm{d}}

\usepackage{xcolor}
\definecolor{turchese}{RGB}{35, 174, 163}

\newtheorem{theorem}{Theorem}

\begin{document}

\begin{frontmatter}

\title{More properties of $(\beta,\gamma)$-Chebyshev functions and points}

\author{Stefano De Marchi}
\ead{demarchi@math.unipd.it}

\author{Giacomo Elefante}
\ead{giacomo.elefante@unipd.it}

\author{Francesco Marchetti}
\ead{francesco.marchetti@unipd.it}

\author{Jean-Zacharie Mariethoz}
\ead{jeanzacharie.mariethoz@studenti.unipd.it}

\address{$^{*}$Dipartimento di Matematica \lq\lq Tullio Levi-Civita\rq\rq, Universit\`a di Padova, Italy}

\begin{abstract}
Recently, $(\beta,\gamma)$-Chebyshev functions, as well as the corresponding zeros, have been introduced as a generalization of classical Chebyshev polynomials of the first kind and related roots. They consist of a family of orthogonal functions on a subset of $[-1,1]$, which indeed satisfies a three-term recurrence formula. In this paper we present further properties, which are proven to comply with various results about classical orthogonal polynomials. In addition, we prove a conjecture concerning the Lebesgue constant's behavior related to the roots of $(\beta,\gamma)$-Chebyshev functions in the corresponding orthogonality interval.
\end{abstract}


\end{frontmatter}


\section{Introduction}
Let $\Omega=[-1,1]$ and $n\in\mathbb{N}$. The well-known Chebyshev polynomials of the first kind $$   T_n(x)=\cos(n\arccos{x}),\;\; n \ge 0,\; x\in\Omega\,, \,$$
are a family of orthogonal polynomials on $\Omega$ with respect to the weight function $w(x)=(1-x^2)^{-1/2}$. $T_n$ is an algebraic polynomial of total degree $n$ with roots
\begin{equation*}
    \mathcal{T}_n=\bigg\{\cos\bigg(\frac{(2j-1)\pi}{n}\bigg)\bigg\}_{j=1,\dots,n},
\end{equation*}
called the {\it Chebyshev points} (of the first kind). These points are pairwise distinct in $\Omega$. Another related set of points in $\Omega$, which includes the extrema of the interval, is the set of \textit{Chebyshev-Lobatto} (CL) points 
\begin{equation*}
    \mathcal{U}_{n+1}=\bigg\{\cos\bigg(\frac{j\pi}{n}\bigg)\bigg\}_{j=0,\dots,n},
\end{equation*}
which consists of the zeros of the polynomial
\begin{equation*}
 \overline{T}_{n+1}(x)=\frac{(1-x^2)}{n}\frac{\mathrm{d}}{\mathrm{d}x}T_n(x),\; x\in\Omega.
\end{equation*}
Chebyshev polynomials and associated points are a classical topic in approximation theory (see, e.g., \cite{Mason02,Rivlin74} for a thorough overview). They represent a valuable choice for constructing stable and effective bases to approximate functions. Indeed, in these tasks Chebyshev and CL points guarantee well conditioning and fast convergence; we refer the interested reader to, e.g., \cite{Cheney98,Rivlin03} for a complete treatment. This interest from a theoretical viewpoint is also motivated by the fact that they have found application in different and various research fields, such as, e.g., numerical quadrature \cite{MasjedJamei05}, solution of differential equations \cite{Sweilam16,Haidvogel79} and group theory \cite{Bircan12}. 

Along with the classical setting, in less and more recent literature many efforts have been made in studying their properties and in expanding their setting by including further \textit{similar} polynomials, generalizations and tools \cite{Caratelli20,Cesarano14,Laine80,Oliveira73,Zhang18}. In particular, in \cite{DeMarchi21a} the new family of $(\beta,\gamma)$-Chebyshev functions (of the first kind) on $\Omega$ was introduced; we recall some important properties of this family in Section \ref{sec:betagamma}, while in later sections we investigate further interesting facts.

In fact, the contribution of this paper can be split in two parts. The first part is devoted to analyzing how the $(\beta,\gamma)$ functions inherit some classical properties of orthogonal polynomials. Specifically, in Section \ref{sec:continued} we provide the connection with the continued fractions, and in Section \ref{sec:generating} we derive the corresponding generating function. Moreover, the Christoffel-Darboux formula and the Sturm-Liouville problem are investigated in Sections \ref{sec:battisti} and \ref{sec:sturm}, respectively. In the second part we focus on the $(\beta,\gamma)$ points. The conditioning of the polynomial interpolation process at such points was investigated in the seminal paper \cite{DeMarchi21a}, where a linear growth with $n$ of the corresponding Lebesgue constant was conjectured for certain configurations of parameters $\beta$ and $\gamma$. In Section \ref{sec:lebesgue}, we provide the proof of such a conjecture, which is crucial to distinguish between well and bad conditioned designs of $(\beta,\gamma)$ nodes in $\Omega$. 

\section{$(\beta,\gamma)$-Chebyshev functions and points}\label{sec:betagamma}

\begin{figure}[H]
  \centering
  \includegraphics[width=0.49\linewidth]{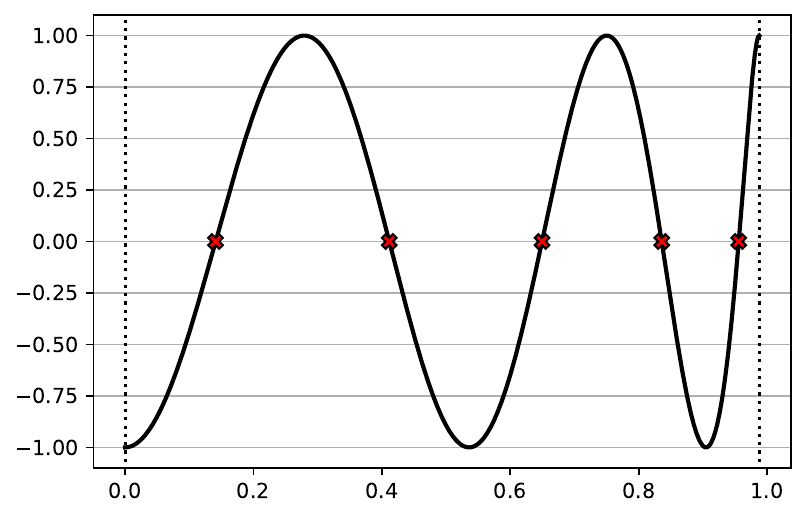} 
\includegraphics[width=0.49\linewidth]{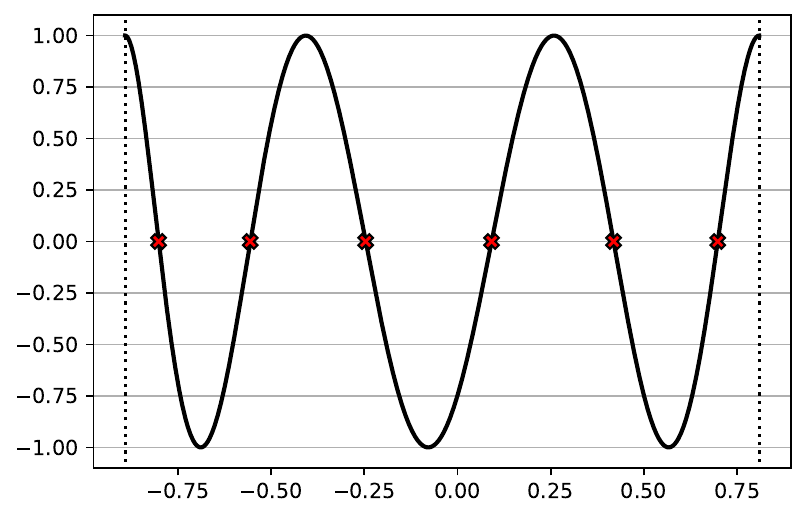}  
\caption{Left: $n=5$, $\beta=1$, $\gamma=0.1$. Right: $n=6$, $\beta=0.3$, $\gamma=0.4$. The points $\mathcal{T}^{\beta,\gamma}_n$ are marked with red crosses.}
\label{fig:0}
\end{figure}

Letting $\beta,\gamma\in[0,2)$, $\beta+\gamma<2$, the family of {\it $(\beta,\gamma)$-Chebyshev functions} (of the first kind) on $\Omega$ consists of functions 
\begin{equation*}\label{eq:betagammacheb}
    T_n^{\beta,\gamma}(x)\coloneqq \cos\bigg(\frac{2n}{2-\beta-\gamma}\bigg(\arccos{x}-\frac{\gamma\pi}{2}\bigg)\bigg),\;x\in\Omega,\; n\in\mathbb{N}.
\end{equation*}
The classical Chebyshev polynomials of the first kind are obtained when $\beta=\gamma=0$, that is $T_n^{0,0}=T_n$, while in general $T_n^{\beta,\gamma}$ is not a polynomial. Covering the same path of the classical framework, we define the \textit{$(\beta,\gamma)$-Chebyshev points} as the zeros of $T_n^{\beta,\gamma}$ in $\Omega_{\beta,\gamma}\coloneqq[-\cos(\beta\pi/2),\cos(\gamma\pi/2)]\subseteq \Omega$, i.e.,
\begin{equation*}
    \mathcal{T}^{\beta,\gamma}_n\coloneqq\bigg\{\cos\bigg(\frac{(2-\beta-\gamma)(2j-1)\pi}{4n}+\frac{\gamma\pi}{2}\bigg)\bigg\}_{j=1,\dots,n}.
\end{equation*}
We show two examples in Figure \ref{fig:0}.

Furthermore, the set of \textit{$(\beta,\gamma)$-Chebyshev Lobatto} points, shortly $(\beta,\gamma)$-CL, are 
\begin{equation*}
    \mathcal{U}^{\beta,\gamma}_{n+1}\coloneqq\bigg\{\cos\bigg(\frac{(2-\beta-\gamma)j\pi}{2n}+\frac{\gamma\pi}{2}\bigg)\bigg\}_{j=0,\dots,n},
\end{equation*}
corresponding to the zeros in $\Omega_{\beta,\gamma}$ of the function 
\begin{equation*}
    \overline{T}_{n+1}^{\beta,\gamma}(x)= \frac{2-\beta-\gamma}{2n} (1-x^2)\frac{d}{dx}T_n^{\beta,\gamma}(x),\;x\in\Omega.
\end{equation*}
Since $\mathcal{T}^{\beta,\gamma}_{n}=\mathcal{U}^{\beta+\frac{2-\beta-\gamma}{2n},\gamma+\frac{2-\beta-\gamma}{2n}}_{n}$ holds true, in the following we will refer directly to $\mathcal{U}^{\beta,\gamma}_{n}$ as $(\beta,\gamma)$-Chebyshev points, with a slight abuse of notation. 

In \cite{DeMarchi21a}, the authors investigated further different aspects of the introduced $(\beta,\gamma)$-Chebyshev functions and points. First, $\mathcal{U}^{\beta,\gamma}_{n}$ can be nicely characterized as a set of \textit{mapped points}. Precisely, considering the Kosloff Tal-Ezer (KTE) map \cite{Adcock16, Kosloff93}
\begin{equation*}\label{kte}
    M_{\alpha}(x)\coloneqq \frac{\sin(\alpha\pi x/2)}{\sin(\alpha\pi /2)},\;\alpha\in]0,1],\;x\in\Omega,
\end{equation*}
and the equispaced points
\begin{equation*}
    \mathcal{E}^{\beta,\gamma}_{n}\coloneqq \bigg\{1-\gamma-\frac{(2-\beta-\gamma)j}{n-1}\bigg\}_{j=0,\dots,n-1},
\end{equation*}
then $\mathcal{U}^{\beta,\gamma}_{n}=M_1\big(\mathcal{E}^{\beta,\gamma}_{n}\big)$. This consideration is crucial for the application of the $(\beta,\gamma)$ setting in the so-called Fake Nodes Approach (FNA), which is a mapped bases scheme that provides well-conditioned reconstruction processes in function approximation with different basis functions \cite{Berrut20,DeMarchi20,DeMarchi21,DeMarchi20b} as well as in numerical quadrature \cite{Cappellazzo22,DeMarchi21b}; we refer to \cite{DeMarchi22} for a recent survey. In \cite{DeMarchi21c}, the set $\mathcal{U}^{\beta,\gamma}_{n}$ is relevant for the analysis of the effectiveness of the {G}ibbs-{R}unge-Avoiding Stable Polynomial Approximation ({GRASPA}) method in the treatment of both Gibbs and Runge's phenomena in the univariate case.

Besides the representation as mapped points, we recall some more properties of the $(\beta,\gamma)$ setting (cf. \cite{DeMarchi21a}).
\begin{itemize}
    \item
    The functions $\{T_n^{\beta,\gamma}\}_{n=0,1,\dots}$ are orthogonal on $\Omega_{\beta,\gamma}$ with respect to the weight function
    \begin{equation*}
    w^{\beta,\gamma}(x)=\frac{2}{(2-\beta-\gamma)\sqrt{1-x^2}},\; x\in \Omega_{\beta,\gamma},
    \end{equation*} 
    precisely
    \begin{equation*}
        \int_{\Omega_{\beta,\gamma}}{T_r^{\beta,\gamma}(x)T_s^{\beta,\gamma}(x)w^{\beta,\gamma}(x)\de x}=
        \begin{dcases} 0 & \textrm{if } r\neq s,\\
        \pi & \textrm{if }r=s=0,\\
        \frac{\pi}{2} & \textrm{if }r=s\neq 0.\end{dcases}
    \end{equation*}
    \item 
    The roles played by $\beta$ and $\gamma$ are symmetric, meaning that $T_n^{\nu_1,\nu_2}$ shows analogous properties as $T_n^{\nu_2,\nu_1}$ as long as $\nu_1,\nu_2$ are \textit{admissible} values for $\beta$ and $\gamma$.
    \item
    For certain values of $\beta$ and $\gamma$ the functions $\{T_n^{\beta,\gamma}\}_{n=0,1,\dots}$ are polynomials, and new results related to the classical setting have been obtained from this analysis; we refer to \cite[\S 2.2]{DeMarchi21a} for further details.
    \item
    Not only Chebyshev ($\beta=\gamma=1/n$) and CL ($\beta=\gamma=0$) points can be defined as $(\beta,\gamma)$-Chebyshev points, but many of their subsets. Indeed, letting $\mathcal{U}_{n+1}=\{u_j\}_{j=0,\dots,n}$, $\kappa_1,\kappa_2\in\mathbb{N}$, $\bs{\kappa}\coloneqq (\kappa_1,\kappa_2),$ and $\beta_{\bs{\kappa}}=\frac{2\kappa_1}{n+\kappa_1+\kappa_2}$, $\gamma_{\bs{\kappa}}=\frac{2\kappa_2}{n+\kappa_1+\kappa_2}$, we obtain
    \begin{equation*}\label{eq:punti_rimanenti}
    \mathcal{U}^{\beta_{\bs{\kappa}},\gamma_{\bs{\kappa}}}_{n} =\mathcal{U}_{n+\kappa_1+\kappa_2}\setminus \{ u_0,\dots,u_{\kappa_2-1}, u_{n+\kappa_2}, \dots, u_{n+\kappa_1+\kappa_2-1}\}.
    \end{equation*}
\end{itemize}

\section{$(\beta,\gamma)$-Chebyshev functions and continued fractions}\label{sec:continued}

It is well-known that any family of orthogonal polynomials $p_1,\dots,p_n,\dots$ satisfies a three-term recurrence formula (see, e.g., \cite[Th. 3.2.1, \S 3.2]{Szego75})
\begin{equation*}
    p_n(x) = (a_n x + b_n)p_{n-1}(x)-c_n p_{n-2}(x),\quad n=2,3,\dots,
\end{equation*}
where $a_n,b_n,c_n$ are constants, $a_n,c_n>0$. Furthermore, if $k_n$ is the coefficient of the $n$ degree monomial in $p_n(x)$, then $a_n=k_n/k_{n-1}$ and $c_n = a_n/a_{n-1}$. For example, Chebyshev polynomials of the first kind satisfy this formula for $n\ge 2$ with $a_n=2$ and $b_n=c_n=1$ (recalling that $T_0(x)=1$ and $T_1(x) = x$).

The functions $\{T_n^{\beta,\gamma}\}_{n=0,1,\dots}$ satisfy a similar recurrence formula in $\Omega_{\beta,\gamma}$
\begin{equation}\label{eq:recurrence}
\begin{split}
    & T_{0}^{\beta,\gamma}(x)=1,\quad T_{1}^{\beta,\gamma}(x)= \cos\bigg(\frac{2}{2-\beta-\gamma}\bigg(\arccos{x}-\frac{\gamma\pi}{2}\bigg)\bigg),\\ & T_{n+1}^{\beta,\gamma}(x)=2T_{1}^{\beta,\gamma}(x)T_{n}^{\beta,\gamma}(x)-T_{n-1}^{\beta,\gamma}(x)
    \end{split}
\end{equation}
Indeed, letting $\theta=\arccos{x}-\frac{\gamma\pi}{2}$, the addition formulae for the cosine easily lead to
\begin{equation*}
    \cos\bigg(\frac{2(n+1)\theta}{2-\beta-\gamma}\bigg)+\cos\bigg(\frac{2(n-1)\theta}{2-\beta-\gamma}\bigg)=2\cos\bigg(\frac{2\theta}{2-\beta-\gamma}\bigg)\cos\bigg(\frac{2n\theta}{2-\beta-\gamma}\bigg).
\end{equation*}
We remark that the only difference with respect to the classical setting relies on the use of the function $T_1^{\beta,\gamma}$ in place of the identity.

Starting from the three terms' recurrence property, our purpose is to characterize the $(\beta,\gamma)$-Chebyshev functions as continued fractions, which are indeed known to be related to orthogonal polynomials.

In the following, we use the space-saving notation 
$$ b_0 + \frac{a_1 \rvert}{\lvert b_1} + \frac{a_2 \rvert}{\lvert b_2} + \frac{a_3 \rvert}{\lvert b_3} + \dots $$
to denote the continued fraction
$$ b_0 + \cfrac{a_1}{b_1 + \cfrac{a_2}{b_2 + \cfrac{a_3}{b_3 + \dots}}}$$
and 
$$ b_0 + \frac{a_1 \rvert}{\lvert b_1} + \frac{a_2 \rvert}{\lvert b_2} + \frac{a_3 \rvert}{\lvert b_3} + \dots + \frac{a_n \rvert}{\lvert b_n}$$
for the $n$-th partial fraction.

Firstly, to deal with the constant factor in the recurrence relation, we introduce the {\it normalized $(\beta,\gamma)$-Chebyshev functions} $\widehat{T}_n^{\beta,\gamma} \coloneqq 2^{1-n} T_n^{\beta,\gamma}$ 
that consequently satisfy the relation
$$ \widehat{T}_n^{\beta,\gamma}(x) = T_1^{\beta,\gamma}(x) \widehat{T}_{n-1}^{\beta,\gamma}(x) - \lambda_n \widehat{T}_{n-2}^{\beta,\gamma}(x),$$
with $\lambda_2 = 1/2$ and $\lambda_n = 1/4, n \geq 3$, $\widehat{T}_{0}^{\beta,\gamma}(x)=1$ and  $\widehat{T}_{-1}^{\beta,\gamma}(x) = 0$.

Then, given $\lambda_1\neq 0$, the continued fraction 
$$ \frac{\lambda_1 \rvert}{\lvert T_1^{\beta,\gamma}(x)} - \frac{\lambda_2 \rvert}{\lvert T_1^{\beta,\gamma}(x)} - \frac{\lambda_3 \rvert}{\lvert T_1^{\beta,\gamma}(x)} - \cdots $$
is such that its $n$-th partial denominator, that is the number $B_n$ such that
$$ \frac{\lambda_1 \rvert}{\lvert T_1^{\beta,\gamma}(x)} - \frac{\lambda_2 \rvert}{\lvert T_1^{\beta,\gamma}(x)} - \frac{\lambda_3 \rvert}{\lvert T_1^{\beta,\gamma}(x)} - \cdots - \frac{\lambda_n \rvert}{\lvert T_1^{\beta,\gamma}(x)} = \frac{A_n}{B_n},$$
satisfies $B_n = \widehat{T}_n^{\beta,\gamma}$. In fact, numerators and denominators of the partial fraction satisfy the recurrence relations (see, e.g., \cite[p. 80]{Chihara11})
\begin{align} \label{RecurAsso}
    B_n &= T_1^{\beta,\gamma}(x) B_{n-1} - \lambda_n B_{n-2}, n \geq 1 \text{ with }B_{-1} = 0, B_0 = 1, \\
    A_n &= T_1^{\beta,\gamma}(x) A_{n-1} - \lambda_n A_{n-2}, n\geq 2  \text{ with }A_{-1} = 1, A_0 = 0 \text{ and } A_1  = \lambda_1. \nonumber
\end{align}
Hence, we can define the function $\widehat{P}_n^{\beta,\gamma} = \lambda^{-1}_1 A_{n+1}$ which is independent of $\lambda_1$ and satisfies the relation
$$ \widehat{P}_n^{\beta,\gamma}(x) = T_1^{\beta,\gamma}(x) \widehat{P}_{n-1}^{\beta,\gamma} - \lambda_{n+1} \widehat{P}_{n-2}^{\beta,\gamma}$$
for $n\geq 1$ and $P_{-1}^{\beta,\gamma}=0,$ $P_0^{\beta,\gamma} = 1$. We refer to these functions as the {\it associated
$(\beta,\gamma)$-Chebyshev functions}. Multiplying the first equation of \eqref{RecurAsso} by $A_{n-1}$ and the second by $B_{n-1}$ we derive the equality (cf. \cite[p. 86]{Chihara11})
\begin{equation}
    \widehat{T}_{n+1}^{\beta,\gamma} \widehat{P}_{n-1}^{\beta,\gamma} -  \widehat{T}_n^{\beta,\gamma} \widehat{P}_n^{\beta,\gamma} = \begin{system}
        - \frac{1}{2} \qquad \qquad\, \text{ for }n=1 \\ \\
        - \frac{1}{2} \big( \frac{1}{4} \big)^{n-2} \quad \text{ for }n\geq 2.
    \end{system}
\end{equation}

\section{The generating function}\label{sec:generating}
Another classical topic is the construction of the generating function associated with a family of functions. 
Suppose that $|u|<1$ and $\phi\in[0,2\pi]$, then (see, e.g., \cite[p. 36]{Rivlin74})
$$ \sum_{n=0}^\infty e^{\iota n \phi}u^n = \sum_{n=0}^\infty (u e^{\iota \phi})^n = \frac{1}{1-u e^{\iota \phi}} $$
By equating the real parts, we get
$$ \sum_{n=0}^\infty  \cos(n\phi)u^n = \frac{1-u \cos \phi}{1+u^2 - 2u \cos \phi}. $$
Similarly to Chebyshev polynomials, we can consider the ordinary generating function related to the $(\beta,\gamma)$ family
$$ F^{\beta,\gamma}(u,x) = \sum_{n=0}^\infty T_n^{\beta,\gamma}(x) u^n\,,$$
and deriving a closed form for $F^{\beta,\gamma}$. Indeed, we can write
$$ T_n^{\beta,\gamma}(x) = \cos\bigg(\frac{2n}{2-\beta-\gamma} \big(\arccos x - \frac{\gamma \pi}{2} \big) \bigg) = \cos (n \phi^{\beta,\gamma}(x)),$$
with $\phi^{\beta,\gamma}(x) \coloneqq \frac{2}{2-\beta-\gamma}(\arccos{x}- \frac{\gamma \pi}{2})$. Then, 
$$ F^{\beta,\gamma}(u,x) = \sum_{n=0}^\infty \cos(n \phi^{\beta,\gamma}(x) ) u^n = \frac{1-u \cos \phi^{\beta,\gamma}(x)}{1+u^2 -2u\cos \phi^{\beta,\gamma}(x)}.$$
Finally, since $\cos \phi^{\beta,\gamma}(x) = T_1^{\beta,\gamma}(x)$, we obtain
$$ F^{\beta,\gamma}(u,x) =  \frac{1-u \, T_1^{\beta,\gamma}(x) }{1+u^2 -2u \,T_1^{\beta,\gamma}(x)}.$$
When $\beta=\gamma=0$, then $T_1^{\beta,\gamma}(x)=x$ and we recover the generating function of the Chebyshev polynomials of the first kind.

\section{The Christoffel-Darboux formula}\label{sec:battisti}
Classical orthogonal polynomials satisfy the Christoffel-Darboux formula. We prove that the same holds for the $(\beta,\gamma)$ functions.

\begin{theorem}
The family of orthogonal functions $T^{\beta,\gamma}_0,T^{\beta,\gamma}_1,\dots,T_n^{\beta,\gamma},\dots$ provides the following Christoffel-Darboux formula
\begin{equation*}
    \sideset{}{'}\sum_{j=0}^n  T^{\beta,\gamma}_j(x)T^{\beta,\gamma}_j(y)=\frac{1}{2}\frac{T^{\beta,\gamma}_{n+1}(x)T^{\beta,\gamma}_n(y)-T^{\beta,\gamma}_n(x)T^{\beta,\gamma}_{n+1}(y)}{T^{\beta,\gamma}_1(x)-T^{\beta,\gamma}_1(y)},
\end{equation*}
with ${\small \displaystyle \sideset{}{'}\sum_{j=0}^n} z_j\coloneqq \frac{1}{2}z_0+z_1+\dots+z_n$.
\end{theorem}

\begin{proof}
We proceed by induction on $n\in\mathbb{N}$. When $n=1$, we need to check the equality
\begin{equation*}    
\frac{1}{2}+T^{\beta,\gamma}_1(x)T^{\beta,\gamma}_1(y)=\frac{1}{2}\frac{T^{\beta,\gamma}_{2}(x)T^{\beta,\gamma}_1(y)-T^{\beta,\gamma}_1(x)T^{\beta,\gamma}_{2}(y)}{T^{\beta,\gamma}_1(x)-T^{\beta,\gamma}_1(y)}.
\end{equation*}
By exploiting the recurrence formula \eqref{eq:recurrence}, we get
\begin{equation*}
\begin{split}
        & \frac{1}{2}\frac{T^{\beta,\gamma}_{2}(x)T^{\beta,\gamma}_1(y)-T^{\beta,\gamma}_1(x)T^{\beta,\gamma}_{2}(y)}{T^{\beta,\gamma}_1(x)-T^{\beta,\gamma}_1(y)}=\\
        & \frac{1}{2}\frac{(2(T^{\beta,\gamma}_{1}(x))^2-1/2)T^{\beta,\gamma}_1(y)-T^{\beta,\gamma}_1(x)(2(T^{\beta,\gamma}_{1}(y))^2-1/2)}{T^{\beta,\gamma}_1(x)-T^{\beta,\gamma}_1(y)}=\\
        & \frac{1}{2}\frac{(1+2T^{\beta,\gamma}_1(x)T^{\beta,\gamma}_1(y))(T^{\beta,\gamma}_1(x)-T^{\beta,\gamma}_1(y))}{T^{\beta,\gamma}_1(x)-T^{\beta,\gamma}_1(y)}=\frac{1}{2}+T^{\beta,\gamma}_1(x)T^{\beta,\gamma}_1(y).
\end{split}
\end{equation*}
Assuming the case $n-1$ holds true, let us prove for $n$. On the left-hand side, we have
\begin{equation*}  
\begin{split}
& \sideset{}{'}\sum_{j=0}^n  T^{\beta,\gamma}_j(x)T^{\beta,\gamma}_j(y)=\\
& \sideset{}{'}\sum_{j=0}^{n-1}  T^{\beta,\gamma}_j(x)T^{\beta,\gamma}_j(y)+T^{\beta,\gamma}_n(x)T^{\beta,\gamma}_n(y)=\\
& \frac{1}{2}\frac{T^{\beta,\gamma}_{n}(x)T^{\beta,\gamma}_{n-1}(y)-T^{\beta,\gamma}_{n-1}(x)T^{\beta,\gamma}_{n}(y)}{T^{\beta,\gamma}_1(x)-T^{\beta,\gamma}_1(y)}+T^{\beta,\gamma}_n(x)T^{\beta,\gamma}_n(y).
\end{split}
\end{equation*}
For the right-hand side, by using the recurrence relation we get
\begin{equation*}
\begin{split}
   & \frac{T^{\beta,\gamma}_{n+1}(x)T^{\beta,\gamma}_n(y)-T^{\beta,\gamma}_n(x)T^{\beta,\gamma}_{n+1}(y)}{2(T^{\beta,\gamma}_1(x)-T^{\beta,\gamma}_1(y))}=\\
   & \frac{(2T^{\beta,\gamma}_{1}(x)T^{\beta,\gamma}_{n}(x)-T^{\beta,\gamma}_{n-1}(x))T^{\beta,\gamma}_n(y)-T^{\beta,\gamma}_n(x)(2T^{\beta,\gamma}_{1}(y)T^{\beta,\gamma}_{n}(y)-T^{\beta,\gamma}_{n-1}(y))}{2(T^{\beta,\gamma}_1(x)-T^{\beta,\gamma}_1(y))}=\\
   & \frac{T^{\beta,\gamma}_{n}(x)T^{\beta,\gamma}_{n-1}(y)-T^{\beta,\gamma}_{n-1}(x)T^{\beta,\gamma}_{n}(y)+T^{\beta,\gamma}_n(x)T^{\beta,\gamma}_n(y)(T^{\beta,\gamma}_1(x)-T^{\beta,\gamma}_1(y))}{2(T^{\beta,\gamma}_1(x)-T^{\beta,\gamma}_1(y))}=\\
   & \frac{T^{\beta,\gamma}_{n}(x)T^{\beta,\gamma}_{n-1}(y)-T^{\beta,\gamma}_{n-1}(x)T^{\beta,\gamma}_{n}(y)}{2(T^{\beta,\gamma}_1(x)-T^{\beta,\gamma}_1(y))}+T^{\beta,\gamma}_n(x)T^{\beta,\gamma}_n(y),
\end{split}
\end{equation*}
which is exactly the left-hand side.
\end{proof}

\section{A Sturm-Liouville problem for $(\beta,\gamma)$ functions}\label{sec:sturm}

Classical orthogonal polynomials can be characterized as solutions of certain differential equations that can be framed into the Sturm-Liouville theory (refer to \cite{AlGwaiz08} for a detailed overview). In the specific case of Chebyshev polynomials of the first kind, $T_n$ is a non-trivial solution of the following boundary-value problem on $\Omega=[-1,1]$
\begin{equation*}
    \frac{\de }{\de x} \Big( \sqrt{1-x^2} \frac{\de y}{\de x} \Big) +\frac{n^2}{\sqrt{1-x^2}}y = 0
\end{equation*}
with {\it separated boundary conditions}
\begin{equation*}
    \begin{dcases}
        a_1 y(-1) + a_2y'(-1) = 0\\
        b_1 y(1) + b_2y'(1) = 0\\        
    \end{dcases}
\end{equation*}
such that $a_1^2+a_2^2>0$ and $b_1^2+b_2^2>0$. Note that $T^\prime_n$ is singular at $x=-1,1$, therefore the boundary equations are satisfied in the limit sense. 

In our framework, we can prove the following
\begin{theorem}
The function $T_n^{\beta,\gamma}$ is a non-trivial solution of the Sturm-Liouville boundary problem on $\Omega_{\beta,\gamma}=[-\cos(\beta\pi/2),\cos(\gamma\pi/2)]\subseteq \Omega$
    \begin{equation} \label{SL-eq}
        \frac{\de }{\de x} \Big( \sqrt{1-x^2} \frac{\de y}{\de x} \Big) + \frac{4n^2}{(2-\beta-\gamma)^2} \frac{1}{\sqrt{1-x^2}}y = 0
    \end{equation}
with separated boundary conditions
    \begin{equation*}
    \begin{dcases}
        a_1 y(-\cos(\beta\pi/2)) + a_2y'(-\cos(\beta\pi/2)) = 0\\
        b_1 y(\cos(\gamma\pi/2)) + b_2y'(\cos(\gamma\pi/2)) = 0\\        
    \end{dcases}
\end{equation*}
where $a_1^2+a_2^2>0$ and $b_1^2+b_2^2>0$.
\end{theorem}

\begin{proof}
We write
\begin{equation*}    
T_n^{\beta,\gamma}(x) = \cos \Big( \frac{2n}{2-\beta-\gamma}\bigg( \arccos(x) - \frac{\gamma \pi}{2}\bigg) \Big) = \cos\Big(\frac{2n}{2-\beta-\gamma} \theta \Big),
\end{equation*}
with $\theta \coloneqq \theta(x) = \arccos(x) - \frac{\gamma \pi}{2}$. By observing that 
\begin{equation} 
\label{Der(b,g)T}
    \frac{\de T_n^{\beta,\gamma} }{\de x} = \frac{\de T_n^{\beta,\gamma} }{\de \theta} \frac{\de \theta }{\de x}, \qquad \frac{\de^2 T_n^{\beta,\gamma} }{\de x^2} = \frac{\de^2 T_n^{\beta,\gamma} }{\de \theta^2} \bigg(\frac{\de \theta }{\de x}\bigg)^2 + \frac{\de T_n^{\beta,\gamma} }{\de \theta} \frac{\de^2 \theta }{\de x^2},
\end{equation}
we get
\begin{align} \label{DerT1}    \frac{\de T_n^{\beta,\gamma} }{\de x} &= \Big(\frac{2n}{2-\beta-\gamma}\Big)\sin\Big(\frac{2n}{2-\beta-\gamma} \theta \Big) \frac{1}{\sqrt{1-x^2}}, \\
    \frac{\de^2 T_n^{\beta,\gamma} }{\de x^2} &= - \frac{4n^2}{(2-\beta-\gamma)^2} \cos\Big(\frac{2n}{2-\beta-\gamma} \theta \Big) \frac{1}{1-x^2} \nonumber \\
    & \qquad +  \Big(\frac{2n}{2-\beta-\gamma}\Big)\sin\Big(\frac{2n}{2-\beta-\gamma} \theta \Big) \frac{x}{(1-x^2)\sqrt{1-x^2}}. \label{DerT2}
\end{align} 

Hence, by substituting $T_n^{\beta,\gamma}$ into the left side of \eqref{SL-eq}, we have
\begin{equation*}
    -\frac{x }{\sqrt{1-x^2}}\frac{\de T_n^{\beta,\gamma}}{\de x} + \sqrt{1-x^2} \frac{\de^2 T_n^{\beta,\gamma} }{\de x^2} + \frac{4n^2}{2-\beta-\gamma} \frac{1}{\sqrt{1-x^2}} T_n^{\beta,\gamma},
\end{equation*}
and from \eqref{DerT1} and \eqref{DerT2} we conclude.
\end{proof}

\section{On the Lebesgue constant of $(\beta,\gamma)$-Chebyshev points for polynomial interpolation}\label{sec:lebesgue}

Let $\mathcal{X}_n\coloneqq \{x_0,\dots,x_n\}$ be a set of distinct points in $\Omega$. The set of Lagrange polynomials $\mathrm{L}\coloneqq\{\ell_0,\dots,\ell_n\}$, where
\begin{equation*}
\ell_i(x)\coloneqq \prod_{\substack{j=0 \\ j\neq i}}^n {\frac{x-x_j}{x_i-x_j}},\; i=0,\dots,n,\;x\in\Omega,
\end{equation*}
allows to construct the {\it Lebesgue function} 
\begin{equation*}
\lambda(\mathcal{X}_n;x)=\sum_{i=0}^n{|\ell_i(x)|},\;x\in\Omega,
\end{equation*}
whose maximum over $\Omega$ is the {\it Lebesgue constant}
\begin{equation*}
    \Lambda(\mathcal{X}_n,\Omega)=\max_{x\in \Omega}\lambda(\mathcal{X}_n;x).
\end{equation*}
As well-known, the Lebesgue constant is an indicator of both the conditioning and the stability of the interpolation process, and it depends on the choice of the interpolation points set $\mathcal{X}_n$. This is why it is important to look for {\it optimal} or nearly-optimal interpolation points, that are points whose Lebesgue constant grows logarithmically in the dimension $n$ of the polynomial space (refer to, e.g., \cite{Cheney98, Rivlin03} for more details). Among other \textit{well-behaved} sets \cite{Ibrahimoglu16}, both Chebyshev and CL points retain the logarithmic growth of the Lebesgue constant \cite{Brutman78,McCabe73}.


Coming back to the $(\beta,\gamma)$ framework, it is reasonable to hypothesize that $(\beta,\gamma)$ points preserve a logarithmic growth of their Lebesgue constant as long as $\beta$ and $\gamma$ are \textit{small enough}, because in this case they can be considered perturbed CL points. In fact, this was proved in \cite[Th. 3]{DeMarchi21a} by exploiting a result in \cite{PiazVian18}. However, as the value of $\beta$ and/or $\gamma$ becomes larger, the behavior of the Lebesgue constant evolves towards an increased rate of growth with $n$. In this sense, it is worthwhile to analyze under which circumstances the growth of the Lebesgue constant related to $(\beta,\gamma)$ points shifts from being logarithmic to be linear. To do so, we consider the setting $\beta=\bar{\beta}_n=2/n$ and $\gamma=0$, which is also interesting since (cf. \eqref{eq:punti_rimanenti}) 
\begin{equation*}
\mathcal{U}^{\bar{\beta}_n,0}_{n}=\mathcal{U}_{n+1}\setminus\{-1\},
\end{equation*}
that is, we are removing the left extremal point from the set of CL points. Now, we prove what was originally conjectured in \cite[Conj. 1]{DeMarchi21a}, i.e., we show that the Lebesgue constant associated to $\mathcal{U}^{\bar{\beta}_n,0}_{n}$ has linear growth and assumes integer values only. Note that in \cite[Th. 4]{DeMarchi21a} we proved the weaker result
$$
  \lambda\big(\mathcal{U}^{\bar{\beta}_n,0}_{n};-1\big)=2n-1.
$$
\begin{theorem}\label{teor_lebesgue}
Let $\bar{\beta}_n=2/n$, $n\in\mathbb{N}_{>0}$. Then, 
\begin{equation*}    \Lambda(\mathcal{U}^{\bar{\beta}_n,0}_{n},\Omega)=2n-1
\end{equation*}
\end{theorem}
\begin{proof}
In light of \cite[Th. 4]{DeMarchi21a}, it is sufficient to prove that the Lebesgue function $\lambda\big(\mathcal{U}^{\bar{\beta}_n,0}_{n};x\big)$ attains its maximum in $x=-1$.
The plan of this technical proof is first to show that the Lebesgue function is monotonically decreasing in $[-1,u_{n-1}]$, where $u_{n-1}$ is the smaller point of the set $\mathcal{U}^{\bar{\beta}_n,0}_{n}$. Then, we conclude by showing that the Lebesgue function is smaller than $2n-1$ in $[u_{n-1},1]$.
    
To this aim, we remind that the Lagrange weight for the $i$-th CL node is (see, e.g., \cite[p. 37]{Trefethen13})
$$ \lambda_i = \prod_{\substack{j=0 \\ j\neq i}}^n \frac{1}{u_i-u_j} = \frac{2^{n-1}}{n} (-1)^i \delta_i$$
where $\delta_i=1$ if $i\neq 0,n$, and it is halved otherwise.
    
Let $x,y\in[-1,u_{n-1}]$ such that $x<y$, then we have
$|x-u_i|>|y-u_i|$ for $i=0,\dots,n-1$. As a consequence for the $i$-th Lagrange polynomial built upon the $\mathcal{U}_{n}^{\bar{\beta}_n , 0}$ nodes, we get 
$$ |\ell_i(x)| = \prod_{\substack{j=0 \\ j\neq i}}^{n-1} \frac{|x-u_j|}{|u_i-u_j|} > \prod_{\substack{j=0 \\ j\neq i}}^{n-1} \frac{|y-u_j|}{|u_i-u_j|} = |\ell_i(y)| .$$
    
Therefore, the Lebesgue function
\begin{align*}
    \lambda(\mathcal{U}_{n}^{\bar{\beta}_n , 0};x) & = \sum_{i=0}^{n-1} |\ell_i(x)| = \sum_{i=0}^{n-1}\prod_{\substack{j=0 \\ j\neq i}}^{n-1} \frac{|x-u_j|}{|u_i-u_j|} \\
    & = \sum_{i=0}^{n-1}\prod_{\substack{j=0 \\ j\neq i}}^{n-1} |x-u_j| |u_i-u_n| \cdot \prod_{\substack{j=0 \\ j\neq i}}^{n} \frac{1}{|u_i-u_j|} \\
    & = \frac{2^{n-1}}{n}  \sum_{i=0}^{n-1}  \Bigg(  \prod_{\substack{j=0 \\ j\neq i}}^{n-1} |x-u_j| (1+u_i)  \delta_i   \Bigg),
\end{align*} 
is strictly decreasing in the interval $[-1,u_{n-1}]$. Now, let us focus on the interval $[u_{n-1},1]$. First, we observe that     
\begin{align*}
      \lambda (\mathcal{U}_{n+1};x) &  = \sum_{i=0}^{n}\prod_{\substack{j=0 \\ j\neq i}}^{n} \frac{|x-u_j|}{|u_i-u_j|} = \sum_{i=0}^{n-1}\prod_{\substack{j=0 \\ j\neq i}}^{n} \frac{|x-u_j|}{|u_i-u_j|} +  \prod_{j=0}^{n-1} \frac{|x-u_j|}{|u_n-u_j|}\\
      &  = \frac{2^{n-1}}{n}\Bigg( \sum_{i=0}^{n-1}  \prod_{\substack{j=0 \\ j\neq i}}^{n-1} |x-u_j|(x+1) \delta_i + \frac{1}{2} \prod_{j=0}^{n-1} |x-u_j| \Bigg).
\end{align*}

Therefore, we get the following
\begin{align*}
      \lambda(\mathcal{U}_{n}^{\bar{\beta}_n , 0};x) & = \lambda(\mathcal{U}_{n}^{\bar{\beta}_n , 0};x) - \lambda (\mathcal{U}_{n+1};x) + \lambda (\mathcal{U}_{n+1};x) \\
      & = \lambda (\mathcal{U}_{n+1};x) + \frac{2^{n-1}}{n} \Bigg( \sum_{i=0}^{n-1} \prod_{\substack{j=0 \\ j\neq i}}^{n-1} |x-u_j| (u_i - x)  \delta_i   -  \frac{1}{2} \prod_{j=0}^{n-1} |x-u_j| \Bigg)\\
      & \le \lambda (\mathcal{U}_{n+1};x) + \frac{2^{n-1}}{n} \Bigg( \sum_{i=0}^{n-1} \prod_{\substack{j=0 \\ j\neq i}}^{n-1} |x-u_j| |x-u_i|  \delta_i   -  \frac{1}{2} \prod_{j=0}^{n-1} |x-u_j| \Bigg)
\end{align*}
Note that when $i=0$ in the summation we get
\begin{equation*}
    \prod_{\substack{j=0 \\ j\neq 0}}^{n-1} |x-u_j| |x-u_0|\cdot\frac{1}{2}   -  \frac{1}{2} \prod_{j=0}^{n-1} |x-u_j| = 0,
\end{equation*}
and we bound
$$ \sum_{i=0}^{n-1} \prod_{\substack{j=0 \\ j\neq i}}^{n-1} |x-u_j| |x-u_i|  \delta_i   -  \frac{1}{2} \prod_{j=0}^{n-1} |x-u_j| \leq \sum_{i=1}^{n-1} \prod_{j=0}^{n-1} |x-u_j|.$$
Therefore, by observing that the terms in the summation do not depend on $i$, we obtain
\begin{equation} \label{BoundLebFun}
    \lambda(\mathcal{U}_{n}^{\bar{\beta}_n , 0};x) \leq \lambda (\mathcal{U}_{n+1};x) + (n-2) \bigg( \frac{2^{n-1}}{n}\prod_{j=0}^{n-1} |x-u_j|\bigg).
\end{equation}
We point out that the term in brackets in \eqref{BoundLebFun} is the double of the absolute value of the $n+1$-th Lagrange polynomial $\ell_n$ constructed at $n+1$ CL nodes. Now, we aim to prove that such value is less or equal than one. By substituting $x = \cos \theta $, and by setting $u_i = \cos \theta_i$ where $\theta_i$ are equidistant points in $[0,\pi]$, we have the following (cf. \cite{Brutman78})
\begin{equation}
    \ell_i (\cos \theta) = \frac{(-1)^i}{2n} \delta_i \sin n\theta \Bigg( \cot \bigg( \frac{\theta + \theta_i}{2} \bigg) + \cot \bigg( \frac{\theta - \theta_i}{2} \bigg) \Bigg).
\end{equation}
For $i=n$, and thus $\theta_i = \pi$, we have
\begin{align*}
    \ell_n (\cos \theta) & = \frac{(-1)^n}{4n} \sin n\theta \Bigg( \cot \bigg( \frac{\theta + \pi}{2} \bigg) + \cot \bigg( \frac{\theta -\pi}{2} \bigg) \Bigg), \\
    & = \frac{(-1)^{n+1}}{2n} \sin n\theta \tan \frac{\theta}{2}.
\end{align*}

By taking the absolute value and by restricting to $x\in [u_{n-1},1]$, i.e., $\theta \in [0,\frac{n-1}{n}\pi]$, we can write
$$ 2|\ell_n(\cos \theta)| = \frac{|\sin n\theta |}{n} \tan \bigg( \frac{\theta}{2} \bigg) \leq \frac{\tan \big( \frac{\theta}{2}\big)}{n} \leq \frac{\tan \big( \frac{n-1}{n}\pi\big)}{n} \leq \frac{1}{n \sin \frac{\pi}{2n}}, $$
since the tangent is strictly decreasing in that interval. Moreover, we observe that the function $g(y) = \frac{1}{y \sin \big( \frac{\pi}{2y}\big)}$ is decreasing for $y\geq 1$, and therefore
\begin{equation} \label{BoundLebFun2}
    2|\ell_n(\cos \theta)|  \leq g(1) = 1.
\end{equation}
Due to \eqref{BoundLebFun} and \eqref{BoundLebFun2} we then have 
\begin{equation}
     \lambda(\mathcal{U}_{n}^{\bar{\beta}_n , 0};x) \leq \lambda (\mathcal{U}_{n+1};x) + n-2.
\end{equation}
Finally, by noticing that the Lebesgue constant of the CL nodes has a logarithmic growth with $n$ increasing, and that
$$ \lambda (\mathcal{U}_{n+1};x) \leq n $$
for $n\ge 1$, we conclude with
$$     \lambda(\mathcal{U}_{n}^{\bar{\beta}_n , 0};x) \leq 2n-2 < 2n -1 .$$
\end{proof}

\begin{figure}[H]
  \centering
  \includegraphics[width=0.49\linewidth]{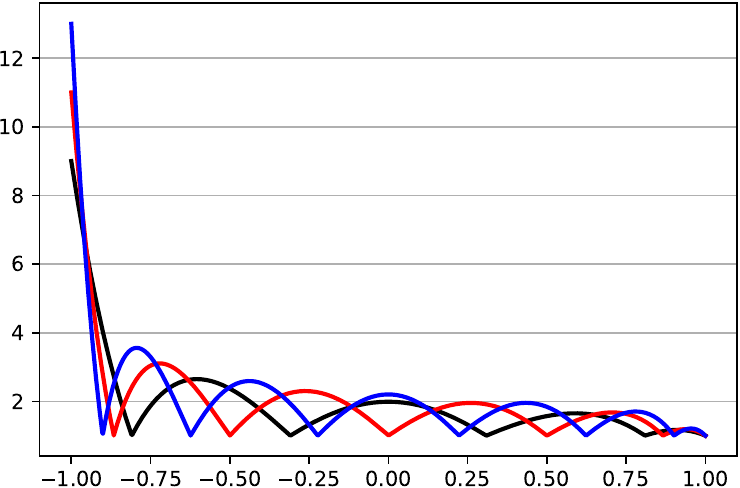}  
\includegraphics[width=0.49\linewidth]{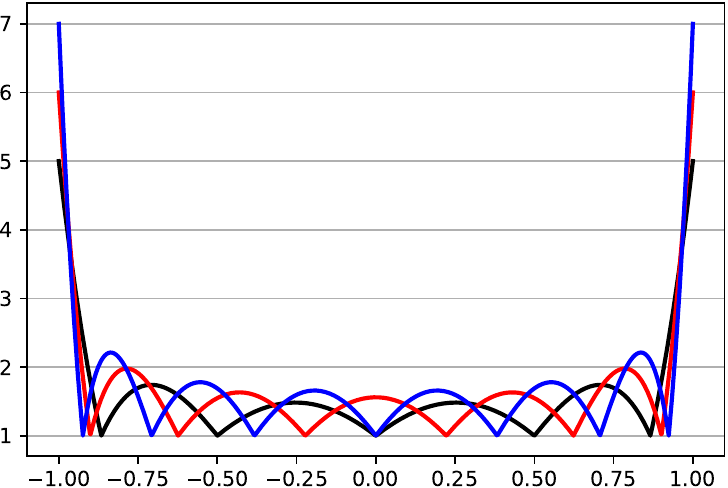}  
\caption{The Lebesgue function associated to $\mathcal{U}^{\bar{\beta}_n,0}_{n}$ (left) and $\mathcal{U}^{\bar{\delta}_n,\bar{\delta}_n}_{n}$ (right). We displayed the cases $n=5$ (black), $n=6$ (red) and $n=7$ (blue).}
\label{fig:1}
\end{figure}

We end this section by observing an interesting analogy with respect to another classical result in approximation theory. As presented, e.g., in \cite{Brutman97}, the Lebesgue constant related to the CL nodes deprived of the extrema of the interval behaves as
\begin{equation*}
\Lambda\big(\mathcal{U}_{n+2}\setminus\{-1,1\},\Omega\big)=n,
\end{equation*}
where $\mathcal{U}_{n+2}\setminus\{-1,1\}=\mathcal{U}^{\bar{\delta}_n,\bar{\delta}_n}_{n}$ with $\bar{\delta}_n=2/(n+1)$. Also in this case, the Lebesgue constant assumes integer values only, and has linear growth with $n$. In Figure \ref{fig:1}, we show some examples of the discussed Lebesgue functions and constants.

\section{Conclusions}

In this work, we showed that $(\beta,\gamma)$-Chebyshev functions satisfy many relevant properties of Chebyshev polynomials of the first kind. After that, we proved a conjecture that concerned the behavior of the Lebesgue constant of $(\beta,\gamma)$-Chebyshev points for certain choices of parameters $\beta$ and $\gamma$ in the polynomial interpolation framework. The obtained results indicate that other classical orthogonal polynomials might be extended to a more general $(\beta,\gamma)$ framework. Furthermore, future work will also focus on the usage of $(\beta,\gamma)$-Chebyshev functions as basis elements in univariate and multivariate approximation schemes.

\section*{Acknowledgments}
This research has been accomplished within GNCS-IN$\delta$AM, Rete ITaliana di Approssimazione (RITA), and the topic group on "Approximation Theory and Applications" of the Italian Mathematical Union (UMI).
The third author also acknowledges the financial support of the Programma Operativo Nazionale (PON) "Ricerca e Innovazione" 2014 - 2020. 

\bibliographystyle{siam}
\bibliography{biblio.bib}
\end{document}